\documentclass[10pt,a4paper,reqno]{amsart}
\usepackage{amsfonts,amsthm,latexsym,amsmath,amssymb,amscd,amsmath,epsf,morefloats}
\usepackage{placeins}
\usepackage{graphicx}
\usepackage{tikz}
\usetikzlibrary{arrows}
\usetikzlibrary{decorations.markings}

\newtheorem{theo+}           {Theorem}
\newtheorem{prop+}           {Proposition}
\newtheorem{coro+}           {Corollary}
\newtheorem{lemm+}           {Lemma}

\theoremstyle{definition}
\newtheorem{defi+}           {Definition}

\theoremstyle{remark}
\newtheorem{rema+}           {Remark}

\newenvironment{theorem}{\begin{theo+}}{\end{theo+}}

\newenvironment{corollary}{\begin{coro+}}{\end{coro+}}
\newenvironment{lemma}{\begin{lemm+}}{\end{lemm+}}

\def \IN{\mathbb N}
\def \IR{\mathbb R}

\numberwithin{equation}{section}

\theoremstyle{remark}

\usepackage[left=1 in,top=1 in,right=1 in,bottom=1 in]{geometry}

\begin{document}

\title[Infinitesimally small spheres and conformally invariant metrics]{Infinitesimally small spheres \\ and \\conformally invariant metrics}

\author[S. Pouliasis]{Stamatis Pouliasis}
\author[A. Yu. Solynin]{Alexander Yu. Solynin}
\address{Texas Tech University-Costa Rica \\ Avenida Escaz\'{u}, Edificio AE205 \\ 
San Jose, Costa Rica, 10203}
\email{{\tt stamatis.pouliasis@ttu.edu}} \email{{\tt
alex.solynin@ttu.edu}}

\date{December 02, 2018}
\thanks{}

\begin{abstract}
The \emph{modulus metric} (also called the \emph{capacity metric})
on a domain $D\subset \mathbb{R}^n$ can be defined as
$\mu_D(x,y)=\inf\{{\mbox{cap}}\,(D,\gamma)\}$, where
${\mbox{cap}}\,(D,\gamma)$ stands for the capacity of the
condenser $(D,\gamma)$ and the infimum is taken over all continua
$\gamma\subset D$ containing the points $x$ and $y$. It was
conjectured by J.~Ferrand, G.~Martin and M. Vuorinen in 1991 that
every isometry in the modulus metric is a conformal mapping. In
this note, we confirm this conjecture and prove  new geometric
properties of surfaces that are spheres in the metric space
$(D,\mu_D)$.
\end{abstract}

\maketitle

\section{Conformal mappings and isometries of the modulus metric} %
A continuous  one-to-one mapping $f:D\to \Omega$ from a domain
$D\subset \mathbb{R}^n$, $n\ge 2$, onto a domain $\Omega\subset
\mathbb{R}^n$ is conformal if it maps smooth curves in $D$ onto
smooth curves in $\Omega$ preserving oriented angles between
intersecting curves. The class of conformal mappings, which is
rich in planar domains (thanks to the Riemann mapping theorem!),
becomes very restrictive in dimensions $\ge 3$. Precisely, by the
classical Liouville's theorem (see, \cite[p. 388]{G}, \cite[p.
19]{V2}, and references therein),  in dimensions $n\ge 3$, every
$C^4$ conformal mapping $f:D\to \Omega$ is a restriction to $D$ of
a M\"{o}bius self-map of $\bar{\mathbb{R}}^n$, where
$\bar{\mathbb{R}}^n$ is the one point compactification of
$\mathbb{R}^n$.

An important characterization of conformal mappings $f:D\to
\Omega$, as well as a characterization of their generalization,
quasiconformal mappings, can be given in terms of their dilatation
$H_f(x)$. Throughout the text we use the following standard
notations. By $|x|$ we denote the Euclidean norm in $\mathbb{R}^n$
and by $S(x,r)=\{y\in \mathbb{R}^n:\, |y-x|=r\}$ and
$B(x,r)=\{y\in \mathbb{R}^n:\,|y-x|<r\}$ we denote, respectively,
the sphere and the open ball in $\mathbb{R}^n$ centered at $x$
with radius $r>0$. We also  use the following shorter notations
$\mathbb{S}=S(\bar{0},1)$ and $\mathbb{B}=B(\bar{0},1)$ for the
unit sphere and unit ball centered at $\bar{0}=(0,\ldots,0)$. Then
the dilatation $H_f(x)$ can be defined as %
\begin{equation}
\label{1.1} H_f(x)=\limsup_{r\to 0}\frac{\max_{y\in
S(x,r)}|f(y)-f(x)|}{\min_{y\in S(x,r)}|f(y)-f(x)|}. %
\end{equation} %
According to a celebrated theorem proved by  Yu.~G.~Reshetnyak in
\cite{R1} (see also Theorem~5.10 in \cite[Chapter II]{R2}) and by
F.~W.~Gehring,
see  Theorem 16 in \cite{G}, a sense preserving homeomorphism $f:D\to \Omega$ is conformal if %
and only if its dilatation $H_f(x)$ is $1$  a.e. on $D$ and
$H_f(x)<\infty$  on $D$. 

Another geometric characterization of conformal mappings can be
given in terms of modules of families of curves. Precisely, a
sense preserving homeomorphism  $f:D\to \Omega$ satisfies
equation~(\ref{1.1}) and therefore it is conformal if and only if
it preserves the modulus of every family $\Gamma$ of curves in
$D$. To define the modulus of a family $\Gamma$ of curves
$\gamma\subset D$, we consider a class ${\mathcal{A}}(\Gamma)$ of
metrics $\rho\ge 0$ admissible for $\Gamma$ in the following
sense: $\rho\in {\mathcal{A}}(\Gamma)$ if and only if $\rho$ is a
non-negative Borel measurable function satisfying  $\int_\gamma
ds\ge 1$ for all locally rectifiable curves $\gamma\in\Gamma$. Now
the modulus
${\mbox{mod}}\,(\Gamma)$ is defined as %
\begin{equation} \label{1.2} %
{\mbox{mod}}\,(\Gamma)=\inf_{\rho\in {\mathcal{A}}(\Gamma)} \int_D
\rho^n\,dm. %
\end{equation}  %

The fact that the modulus ${\mbox{mod}}\,(\Gamma)$ defined by
(\ref{1.2}) is conformally invariant is classical, see \cite{G}
and \cite[p. 54]{V2}. But verification of invariance of modules of
all families of curves under a mapping $f:D\to \Omega$ is
impractical and there is no need to verify invariance of modules
of every family of curves. It was shown by Gehring \cite{G} that a
mapping will be conformal if it preserves modules of curves
connecting boundary components of the so-called ring domains.
Later, J.~Ferrand, G.~Martin and M.~Vuorinen suggested in
\cite{FMV} that it might be sufficient to verify invariance of
modules for some other specific families of curves. On this way,
these authors studied in \cite{FMV} a conformal invariant $\mu_D$,
which can be defined as follows. If $D$ is a domain in
$\mathbb{R}^n$ and $x,y\in D$, then %
\begin{equation}  \label{1.3} %
\mu_D(x,y)=\inf\{{\mbox{mod}}\,(\Gamma(D,l)):\, l\in C_{xy}\},
\end{equation} %
where $C_{xy}$ is the family of all Jordan arcs $l$ joining $x$ to
$y$ in $D$ and $\Gamma(D,l)$ is the family of all curves $\gamma$
(not necessarily Jordan) in $D$ joining $l$ and $\partial D$. It
was mentioned in \cite[p. 103]{V2} that the function $\mu_D(x,y)$
defines a metric on $D$, which is called the \emph{modulus
metric}, if and only if $\partial D$ is of positive conformal
capacity. Thus, a domain $D$ supplied with the modulus metric
$\mu_D$ becomes a metric space $(D,\mu_D)$.

The modulus metric is conformally invariant and quasi-invariant
under quasiconformal mappings which makes it very useful in the
theory of quasiconformal mappings; see, for instance, \cite{V1},
\cite{V2}. It was conjectured by Ferrand, Martin and Vuorinen (see
\cite[p. 195]{FMV}) that every mapping $f:D\to \Omega$, which is
an \emph{isometry} with respect to the metrics $\mu_D$ and
$\mu_\Omega$, is conformal. These authors have shown in \cite{FMV}
that this is indeed the case when $D$ is a ball in $\mathbb{R}^n$.
An essential progress towards the solution of this conjecture was
made recently in \cite{BP2}, where the authors proved that $f$ is
conformal if $n=2$, thus settling the conjecture in this case, and
that $f$ is quasiconformal in the case $n\ge 3$. The main goal of
this paper
is to prove the following. 
\begin{theorem}  \label{Theorem-1} %
Let $D$ and $\Omega$ be domains in $\mathbb{R}^n$, $n\ge 2$, such
that $\partial D$ and $\partial \Omega$ have positive conformal
capacities. Suppose that $f$ is an isometry of the metric space
$(D,\mu_D)$ onto the metric space $(\Omega,\mu_\Omega)$. Then $f$
is a conformal mapping.
\end{theorem} %

Thus, Theorem~1 proves Ferrand-Martin-Vuorinen conjecture in all
dimensions. The method used to prove this result is purely
geometrical. It is based on application of polarization
transformation in the spirit of papers \cite{S1} and \cite{S4},
where polarization was used to solve P\'{o}lya-Szeg\"{o} problem
on continuous symmetrization.

One more goal of this work is to study geometric  properties of
the $\mu_D$-spheres, that are level surfaces with respect to the
modulus metric, defined by $S_{\mu_D}(x,r)=\{y\in D:\,
\mu_D(x,y)=r\}$ and the $\mu_D$-balls defined by
$B_{\mu_D}(x,r)=\{y\in D:\,\mu_D(x,y)<r\}$. In this direction we
prove the following result.
\begin{theorem}  \label{Theorem-2} %
Let $D$  be a domain in $\mathbb{R}^n$, $n\ge 2$, such that
$\partial D$  has positive conformal capacity and let $x_0\in D$.
Then there is $r_0>0$ such that for all $0<r<r_0$, the
$\mu_D$-sphere $S_{\mu_D}(x_0,r)$ is a topological sphere in
$\mathbb{R}^n$ which satisfies the interior and exterior cone
conditions and the ball $B_{\mu_D}(x_0,r)$  is starlike with
respect to $x_0$.
\end{theorem} %

The rest of the paper is organized as follows. Section~2 contains
necessary background from Potential Theory. In Section~3, a
geometric transformation called \emph{polarization} will be used
to establish some properties of the modulus metric and geometric
properties of $\mu_D$-balls stated in Theorem~2. In Section~4, we
use Reshetnyak's characterization of  conformality by an
invariance property of collections of infinitesimally small
spheres and the lemma on the three spheres from elementary
geometry to prove our Theorem~1.  Finally, in Section~5, we
discuss some related open problems.


\section{Condenser capacity and modulus metrics}
A condenser in $\overline{\mathbb{R}}^n$, $n\geq 2$, is a pair
$(D,K)$, where $D$ is a domain in $\overline{\mathbb{R}}^n$ and
$K$ is a non-empty compact subset of $D$. For every condenser
$(D,K)$, we denote by $H(D,K)$ the class of functions admissible
for $(D,K)$; i.e. $H(D,K)$ consists of all $C^1$ functions with
compact support in $D$ satisfying $u(x)\geq 1$ for $x\in K$. If
$D\subset \IR^n$, the conformal capacity of the condenser $(D,K)$
is defined by
\begin{equation}\label{2.1}
{\rm cap}(D,K)=\inf_u\int_{D\setminus K} |\nabla u|^n\;dm,
\end{equation}
where $m$ is the $n$-dimensional Lebesgue measure and the infimum
is taken over the class $H(D,K)$. If $D$ contains the point
$\infty$, ${\rm cap}(D,K)$ is defined by means of an auxiliary
M\"{o}bius transformation. Since we do not use other capacities in
this paper, everywhere below we  use a shorter term ``capacity''
instead of ``conformal capacity''.

It is instrumental for us, that the capacity of a condenser
$(D,K)$ is conformally invariant, see \cite{HKM}, and, due to
Ziemer's theorem \cite[Theorem~3.8]{Z} (see also
\cite[Theorem~1]{G} and Proposition~10.2 in \cite{Rick}), it
coincides with the modulus of the family $\Gamma(D,K)$ of
curves $\gamma$ joining $K$ and $\partial D$ in $D$; i.e.,  %
$$ 
{\rm cap}(D,K)={\mbox{mod}}\,(\Gamma(D,K)).
$$ 
Therefore, the modulus metric $\mu_D$ can be alternatively defined as %
\begin{equation}\label{2.2}
\mu_D(x,y)=\inf\{{\rm cap}(D,\gamma):\gamma\in C_{xy}\}.
\end{equation}

Thus, the modulus metric is, in a certain sense, the
``\emph{capacity metric}''. The latter definition has some
advantages when studying properties of the metric. Next, we
introduce necessary terminology and recall several known or
``semi-known''
properties of capacities. It would be convenient to list 
these properties as a series of lemmas.

First, we recall that a compact set $E\subset
\overline{\mathbb{R}}^n$ is said to be of \emph{zero capacity} if
there exists a domain $D$ with $E\subset D$ such that
${\mbox{cap}}\,(D,E)=0$. Otherwise, $E$ is said to be a set of
positive capacity (or a set of positive conformal capacity as we
state it in Theorem~1).

It is well known that the infimum in the definition of the
capacity  can be taken, with the same result, over different
classes of functions. Indeed, we recall first that every closed
set $E\subset\IR^{n}$ may contain points which are regular for the
Dirichlet problem for the $n$-Laplacian and points which are
irregular for this problem, see, for instance, Section~9.5 in
\cite{HKM} for the definition and the properties of irregular
boundary points. One but not both of these sets can be empty. Let
$I(E)$ denote the set of irregular points of $E$ for the problem
under consideration.
Let $C(D)$ be the space of continuous functions in $D$, let
$C^{\infty}(D)$ be the space of infinitely differentiable
functions in $D$ and let $C_{0}^{\infty}(D)$ be the subspace of
$C^{\infty}(D)$ of functions with compact support in $D$. Let
$H^{1,n}(D)$ be the completion of $C^{\infty}(D)$ and let
$H_{0}^{1,n}(D)$ be the completion of $C_{0}^{\infty}(D)$ with
respect to the norm %
$$  %
||v||_{1,n}=\Big(\int_{D}|v|^{n}\,dm\Big)^{1/n}+
\Big(\int_{D}|\nabla v|^{n}\,dm\Big)^{1/n},\qquad v\in
C^{\infty}(D).  %
$$  %
 Also, let $W^{1,n}(D)=H^{1,n}(D)\cap C(D)$ and
$W_{0}^{1,n}(D)=H_{0}^{1,n}(D)\cap C(D)$.

\begin{lemma}[see {\cite[pp. 28-29]{HKM}}] \label{Lemma-1} %
Let $(D,K)$ be a condenser in $\IR^{n}$. Then
\begin{equation} \label{2.3} %
{\rm cap}(D,K)=\inf_{u\in H_{1}(D,K)}\int_{D\setminus K} |\nabla u|^n\;dm,%
\end{equation}   %
where $H_{1}(D,K)$ is the family of  functions $u$ in the space
$W_0^{1,n}(D)$ satisfying $u\ge 1$ on $K$.
\end{lemma}

We can restrict the classes of functions over which the infima in
(\ref{2.1}) and (\ref{2.3}) are taken to the subclasses of
so-called \emph{monotone functions}. A continuous function $f$ on
a domain $D\subset\IR^{n}$ is called monotone on $D$ if for any
relatively compact domain $\Omega$ in
$D$,  %
\begin{equation}  \label{2.4} %
\sup_{x\in\partial \Omega}f(x)=\sup_{x\in \Omega}f(x)\qquad\mbox{and}\qquad
\inf_{x\in\partial \Omega}f(x)=\inf_{x\in \Omega}f(x).  %
\end{equation} %

Given a domain $D$ and a compact set $K\subset D$, by $H^*(D,K)$
we denote the family of  monotone on $D\setminus K$ functions $u$
in the space $W_0^{1,n}(D)$ satisfying $u\ge  1$ on $K$.

\begin{lemma}[cf. {\cite[p. 54]{Rick}}]\label{Lemma-2} %
Let $(D,K)$ be a condenser in $\IR^{n}$. Then
$$ 
{\rm cap}(D,K)=\inf_{u\in H^*(D,K)}\int_{D\setminus K} |\nabla u|^n\;dm.  %
$$ 
\end{lemma}

It is well known that every condenser $(D,K)$ has a unique
potential function $\omega_{D,K}$. %
In our next lemma, we summarize some well-known  properties of
this function.  %

\begin{lemma}[see, {\cite[pp. 194,211,212]{HKM}},{\cite[p. 104]{Y}}]\label{Lemma-3} %
For every condenser $(D,K)$ with positive capacity there is a
function $\omega_{D,K}$, called the potential function, which
minimizes the integral in (\ref{2.3}). The potential function
$\omega_{D,K}$ is a solution to the $n$-Laplace equation in
$D\setminus K$ and satisfies the following boundary conditions %

$$  %
\omega_{D,K}(x)\to 1 \quad {\mbox{as $x\in D\setminus K$
approaches regular points of $K$}}
$$ %
and %
$$  %
\omega_{D,K}(x)\to 0 \quad {\mbox{as $x\in D\setminus K$
approaches regular points of $\partial D$.}}
$$ %

Furthermore, if $D$ is bounded, then $\omega_{D,K}\in
H_0^{1,n}(D)$.
\end{lemma} %

The potential function $\omega_{D,K}$ possessing properties
described in Lemma~3 is, in fact, the \emph{$n$-harmonic measure}
of the compact set $K$ with respect to the domain $D\setminus K$.
In particular, $\omega_{D,K}$ is monotone on $D\setminus K$.
For further properties of the $n$-harmonic measure the reader may
consult \cite[Chapter~11]{HKM}. The following convergence lemma is a useful
tool often used to prove existence of condensers with special
properties.

\begin{lemma}[{\cite[Theorem 2.2]{HKM}}]  \label{Lemma-4} %
Let $(D_k,K_k)$, $k\in \mathbb{N}$, be a sequence of condensers in
$\overline{\mathbb{R}}^n$ such that $D_k\subset D_{k+1}$,
$K_{k+1}\subset K_k$ for all $k\in \mathbb{N}$, $\cup_{k=1}^\infty
D_k=D$, $\cap_{k=1}^\infty K_k=K$, and $(D,K)$ is a condenser of
positive capacity.  Then %
$$ %
{\mbox{cap}}\,(D_k,K_k)\to {\mbox{cap}}\,(D,K), \quad {\mbox{as
$k\to \infty$.}}
$$ %
\end{lemma}  %

We will need the following monotonicity property of the capacity.

\begin{lemma}\label{Lemma-5} 
Let $(D_1,K_1)$ and $(D_2,K_2)$ be  condensers in $\IR^{n}$ such
that $D_2\subset D_1$ and $K_1\subset K_2$. Then
\begin{equation}  \label{2.5} %
{\rm cap}(D_1,K_1)\le {\rm cap}(D_2,K_2).  %
\end{equation}   %
Furthermore, if $I(K_1)=\emptyset$ and $K_2\setminus K_1$ contains
a compact set $K_3$  such that $I(K_3)=\emptyset$ and if, in
addition, $D_1\setminus K_1$ is connected and contains $K_3$ and
$D_1\setminus K_3$ is connected and contains $K_1$ then
(\ref{2.5}) holds with the sign of strict inequality.%
\end{lemma} %

\begin{proof} %
The non-strict inequality (\ref{2.5}) is well known, see, for
instance, \cite[Theorem~2.2]{HKM}, and follows immediately from
the definition (2.1). Thus, we have to prove only the statement
about the cases of equality. Also, since the capacity of a
condenser defined by (\ref{2.1}) is conformally invariant we may
assume without loss of generality that $\infty \not\in D_1$.

For $\delta>0$, let $K_3(\delta)=\{x\in \mathbb{R}^n:\,
{\mbox{dist}}\,(x,K_3)<\delta\}$. If
$\delta<{\mbox{dist}}\,(K_3,\partial(D_1\setminus K_1))$, then
$\overline{K_3(\delta)}$ and $\partial (K_3(\delta))$ are compact
subsets of $D_1\setminus K_1$. Let $\omega$ be the potential
function of the condenser
$(D_1,K_1\cup K_3)$ and let %
$$ %
t_m=\max \{\omega(x):\,x\in
\partial (K_3(\delta))\}.  %
$$ %
 It follows from the maximum principle for solutions of the
$n$-Laplace equation (see, for instance, \cite[p. 115]{HKM}) that
$0<t_m<1$. Let $t_m<t_\delta<1$. Our assumptions imply  that
$\omega$ is continuous on the
set $\overline{K_3(\delta)}$ and therefore the set %
$V=\{x\in \overline{K_3(\delta)}:\,\omega(x)>t_\delta\}$ is open
in $\mathbb{R}^n$ and its complement
$E=\overline{\mathbb{R}}^n\setminus V$ is compact in
$\overline{\mathbb{R}}^n$.

Since $\Omega=\overline{\mathbb{R}}^n\setminus K_3$ is open and
connected and $E\subset \Omega$, the pair $(\Omega,E)$ is a
condenser with the
potential function $\omega_E$ given by %
$$ %
\omega_E(x)=\left\{ %
{\begin{array}{cl}\frac{1-\omega(x)}{1-t_\delta}&\quad {\mbox{if
$x\in \Omega\setminus E$}}\\ %
1 &\quad {\mbox{if $x\in E$.}} %
\end{array} }\right.%
$$ %


Let $\varepsilon>0$ be as small as we will need it later. It
follows from the convergence Lemma~4 that there exists a domain
$D_\varepsilon$ such that $K_2\subset D_\varepsilon\subset
\overline{D_\varepsilon}\subset D_1$, the boundary $\partial
D_\varepsilon$
is regular for the $n$-Laplace equation, and %
\begin{equation}  \label{2.6} %
{\mbox{cap}}(D_\varepsilon,K_1\cup
K_3)<{\mbox{cap}}\,(D_1,K_1\cup K_3)+\varepsilon. %
\end{equation}  %
Let $\omega_\varepsilon$ denote the potential function of the
condenser  $(D_\varepsilon,K_1\cup K_3)$. It follows from the
maximum principle for
solutions of the $n$-Laplace equation that %
$$ 
\omega_\varepsilon(x)\le \omega(x) \quad {\mbox{for all $x\in
D_\varepsilon$.}}
$$ 
The latter equation implies that the set
$V_{\varepsilon,\delta}=\{x\in
\overline{K_3(\delta)}:\,\omega_\varepsilon(x)>t_\delta\}$ is an
open subset of $V$, the set
$E_{\varepsilon,\delta}=\overline{\mathbb{R}}^n\setminus
V_{\varepsilon,\delta}$ is compact in $\overline{\mathbb{R}}^n$,
and $E_{\varepsilon,\delta}\supset E$. Furthermore, the pair $(\Omega,E_{\varepsilon,\delta})$
can be considered as a condenser and the function %
$$ %
\omega_{\varepsilon,\delta}(x)=\left\{ %
{\begin{array}{cl}\frac{1-\omega_\varepsilon
(x)}{1-t_\delta}&\quad {\mbox{if
$x\in \Omega\setminus E_{\varepsilon,\delta}$}}\\ %
1 &\quad {\mbox{if $x\in E_{\varepsilon,\delta}$}} %
\end{array} }\right.%
$$ %
is the potential function of  $(\Omega,E_{\varepsilon,\delta})$.
This yields the following inequality: %
\begin{equation} \label{2.7}
\int_{V_{\varepsilon,\delta}\setminus K_3}|\nabla
\omega_\varepsilon|^n\,dm=(1-t_\delta)^n
{\mbox{cap}}\,(\Omega,E_{\varepsilon,\delta})\ge (1-t_\delta)^n
{\mbox{cap}}\,(\Omega,E)=\int_{\Omega\setminus E}|\nabla
\omega|^n\,dm.
\end{equation}  %

Notice that the function %
$$ %
u(x)=\left\{ %
{\begin{array}{cl}\omega_\varepsilon (x)&\quad {\mbox{if
$x\in D_\varepsilon\setminus V_{\varepsilon,\delta}$}}\\ %
t_\delta &\quad {\mbox{if $x\in V_{\varepsilon,\delta}$}} %
\end{array} }\right.%
$$ %
is admissible for the condenser $(D_\varepsilon,K_1)$. This together with relations (\ref{2.6}) and (\ref{2.7}) implies  %

\begin{alignat}{10} \label{2.8} %
{\mbox{cap}}\,(D_1,K_1)\le {\mbox{cap}}\,(D_\varepsilon,K_1)\le
\int_{D_\varepsilon\setminus (K_1\cup K_3)}|\nabla
\omega_\varepsilon|^n\,dm-\int_{V_{\varepsilon,\delta}\setminus
K_3}|\nabla \omega_\varepsilon|^n\,dm   \\ %
{ } ={\mbox{cap}}\,(D_\varepsilon,K_1\cup
K_3)-(1-t_\delta)^n{\mbox{cap}}\,(\Omega,E_{\varepsilon,\delta}) \
\ \ \ \ \ \ \ \ \ \ \ \ \ \ \ \ \ \ \ \ \ \ \   \nonumber
\\ \le {\mbox{cap}}\,(D_1,K_1\cup K_3)+\varepsilon-(1-t_\delta)^n
{\mbox{cap}}\,(\Omega,E).  \ \ \ \ \ \ \ \ \ \ \ \ \ \ \ \ \ \ \ \
\     \nonumber
\end{alignat} %

Finally, assuming that  $\varepsilon<(1-t_\delta)^n
{\mbox{cap}}\,(\Omega, E)$ and using the non-strict monotonicity
property of capacity of condensers, we conclude from
(\ref{2.8}) that %
$$  %
{\mbox{cap}}\,(D_1,K_1)<{\mbox{cap}}\,(D_1,K_1\cup K_3)\le
{\mbox{cap}}\,(D_1,K_2)\le {\mbox{cap}}\,(D_2,K_2), %
$$ %
which is the required strict monotonicity property.  \end{proof}  %


\smallskip




We note here that the value of the modulus metric $\mu_D(x,y)$
does not change if we replace $C_{xy}$ with the family $K_{xy}$ of
all continua (connected compact sets) in $D$ containing $x,y$.
Precisely, we have the following result.

\begin{lemma}[see {\cite[p. 191]{FMV}}]\label{Lemma-6}  %
The modulus metric $\mu_D$ can be alternatively defined as %
\begin{equation}\label{2.9}
\mu_D(x,y)=\inf\{{\rm cap}(D,K):K\in K_{xy}\}.
\end{equation}  %
\end{lemma}      %
 \begin{proof}
Let $\nu_D(x,y)$ denote the infimum in (\ref{2.9}). Since
$C_{xy}\subset K_{xy}$ it is immediate from (\ref{2.2}) and
(\ref{2.9}) that  $\nu_D(x,y)\le \mu_D(x,y)$.

To prove the reverse inequality we consider a sequence of continua
$K_k$, $k\in \mathbb{N}$,    in $K_{xy}$ such that
${\mbox{cap}}\,(D,K_k)\to \nu_D(x,y)$ as $k\to \infty$ and a
sequence of $\varepsilon_k>0$ such that $\varepsilon_k\to 0$ as
$k\to \infty$. It follows from the convergence Lemma~4 that for
every $k\in \mathbb{N}$ there exists $\delta_k>0$ such that the
set $K(\delta_k)=\{x\in R^n:\, {\mbox{dist}}\,(x,K_k)\le
\delta_k\}$ is a compact subset of $D$ and %
\begin{equation}  \label{2.10} %
{\mbox{cap}}\,(D,K(\delta_k))\le
{\mbox{cap}}\,(D,K_k)+\varepsilon_k. %
\end{equation} %
Furthermore, the interior of $K(\delta_k)$ is a non-empty
connected open set containing points $x$ and $y$. Therefore, for
every $k\in \mathbb{N}$,  there exists a Jordan arc
$\gamma_k\subset K(\delta_k)$ (one may assume that it is analytic
if necessary) joining points $x$ and $y$. Now, by (\ref{2.5}) and
(\ref{2.10}), %
$$ %
\mu_D(x,y)\le {\mbox{cap}}\,(D,\gamma_k)\le
{\mbox{cap}}\,(D,K_k)+\varepsilon_k.
$$ %
Taking the limit in the last inequality we obtain $\mu_D(x,y)\le
\nu_D(x,y)$, which combined with the reverse inequality
mentioned above gives (\ref{2.9}) %
 \end{proof} %



An advantage of the definition of $\mu_D(x,y)$ given by
(\ref{2.9}) is that it is easier to establish existence of a
continuum minimizing the capacity in the right-hand side of
(\ref{2.9}) than to prove that this extremal continuum is a Jordan
arc. For instance, Lemma~6 below guarantees, in most cases,
existence of a continuum extremal for problem (\ref{2.2}) but does
not provide enough information to conclude that this continuum  is
a Jordan arc extremal for problem (\ref{1.3}). Similar existence
results are known for some other problems (see \cite{F}) but for
the problem under consideration it was not recorded in the
literature available for us. Thus, we provide its proof here.

\begin{lemma}\label{Lemma-7}
Let $D$ be a domain in $\IR^{n}$ such that $\partial D$ has
positive capacity, let $E\subset D$ be a connected compact set and
let $x,y\in E$. Suppose that there is a sequence of Jordan arcs
$\gamma_k\in C_{xy}$ such that $\gamma_k\subset E$ for all $k\in
\mathbb{N}$ and ${\mbox{cap}}\,(D,\gamma_k)\to \mu_D(x,y)$ as
$k\to \infty$.
Then there exists a continuum $\beta\subset E$ containing $x$ and
$y$ such that
\begin{equation}\label{2.11}
    \mu_D(x,y)={\rm cap}(D,\beta).
\end{equation}
\end{lemma} %

\begin{proof} %
Let $\omega_{k}=\omega_{D,\gamma_k}$ be the potential function of
the condenser $(D,\gamma_k)$ that is also the $n$-harmonic measure
of $\gamma_{k}$ with respect to  $D\setminus\gamma_{k}$. From the
maximum and minimum principles for $n$-harmonic functions
and Corollary~2.5 in \cite{F} 
we conclude that for every $k\in\IN$ the function $\omega_{k}$ is
monotone (in the sense of definition (\ref{2.4}))  on the domain
$D\setminus\{x\}$.
Using this fact, Proposition~1.6 in \cite{F} and 
passing to a subsequence if necessary, we may assume that the
sequence of functions $\omega_{k}$ converges locally uniformly on
$D\setminus\{x\}$ to a continuous function $\omega$ which has
generalized
partial derivatives on $D\setminus\{x\}$, satisfying %
\begin{equation}\label{2.12}
\int_D |\nabla \omega|^n\;dm\leq\lim_{k\to\infty}\int_{D\setminus
\gamma_{k}} |\nabla \omega_{k}|^n\;dm=\lim_{k\to\infty}{\rm
cap}(D,\gamma_{k})=\mu_D(x,y).
\end{equation}
Let $\omega_E$ be the $n$-harmonic measure of $E$ with respect to
$D\setminus E$. Since $\gamma_k\subset E$,  the Carleman's
principle for $n$-harmonic measures (see Theorem~11.3 in
\cite{HKM}) implies that $0\leq\omega_{k}\leq \omega_E$ on
$D\setminus E$, for every $k\in\IN$. Therefore, letting
$k\to\infty$, we conclude that $0\leq\omega\leq \omega_E$ on
$D\setminus E$. For every regular boundary point $\zeta\in\partial
D$, we have
$$0\le \liminf_{z\to \zeta}\omega(z)\le \limsup_{z\to \zeta}\omega(z)
\le  \lim_{z\to \zeta} \omega_E(z)=0.$$
Therefore, %
$$ 
\omega(z)\to 0, \mbox{ as }z\to \partial D\setminus I(\partial D). %
$$ 

Let $\beta=\omega^{-1}(1)\subset E$. Suppose that $\beta$ is not
connected. Let $S\subset D\setminus\beta$ be a topological sphere
such that both connected components $\Omega_{1}$ and $\Omega_{2}$ of $\IR^{n}\setminus S$
intersect $\beta$. Note that, since $S$ is compact,
 \begin{equation} \label{2.13}
\max_{z\in S}\omega(z)<1.
\end{equation}
Suppose that there is a subsequence $\gamma_{k_{m}}$
of $\gamma_{k}$ such that $\gamma_{k_{m}}\cap S=\emptyset$,
for all $m\in\IN$. We may assume, passing to a subsequence of 
$\gamma_{k_{m}}$ if needed, that $\gamma_{k_{m}}$, $m\in\IN$, lies in the same
component of $\IR^{n}\setminus S$, say $\Omega_{1}$. Then $\omega_{k_{m}}$
are $n$-harmonic functions on $\Omega_{2}$ which converge locally uniformly
to $\omega$ on $\Omega_{2}$. From Theorem 6.13 in
\cite[p. 117]{HKM}, $\omega$ is $n$-harmonic on $\Omega_{2}$.
Let $w_{0}\in\Omega_{2}\cap\beta\neq\emptyset$. Then
$$\omega(w_{0})=1\geq\omega(z),$$
for all $z\in\Omega_{2}$. From the maximum principle for
$n$-harmonic functions (see, for instance, \cite[p. 115]{HKM}),
$\omega=1$ on $\Omega_{2}$. Since $\omega$ is continuous on
$\Omega_{2}\cup S$, we obtain that $\omega=1$ on $S$, which
contradicts (\ref{2.13}).
We conclude that
there exists $p\in\IN$ such that $\gamma_{k}\cap
S\neq\varnothing$ for every $k\geq p$. Let $z_{k}\in\gamma_{k}\cap
S$, $k\geq p$. Since S is compact we may assume that $z_{k}\to
z_{0}\in S$. Since $\omega_{k}\to\omega$ uniformly on $S$,
$$\omega(z_{0})=\lim_{k\to\infty}\omega_{k}(z_{k})=1,$$
contradicting (\ref{2.13}). Therefore $\beta$ is connected. Since
obviously $\beta$ is closed, $\beta\in K_{xy}$.

Let $z\in D\setminus\beta$. Since $\omega(z)<1$ and $\omega_{k}=1$
on $\gamma_{k}$, there exists $\varepsilon>0$ and $k_{0}\in\IN$
such that $B(z,\varepsilon)\cap\gamma_{k}=\emptyset$ for every
$k\geq k_{0}$. Therefore, for every $k\geq k_{0}$, $\omega_{k}$ is
$n$-harmonic on $B(z,\varepsilon)$ and $\omega_{k}\to\omega$
locally uniformly on $B(z,\varepsilon)$. From Theorem 6.13 in
\cite[p. 117]{HKM}, $\omega$ is $n$-harmonic on
$B(z,\varepsilon)$. Since $z\in D\setminus\beta$ was arbitrary,
$\omega$ is $n$-harmonic on $D\setminus\beta$. We conclude that
$\omega$ is $n$-harmonic on $D\setminus\beta$ with boundary values
$1$ on $\beta$ and $0$ on every regular boundary point of $D$.
From Theorem 11.1(c) \cite[p. 209]{HKM}) we get that $\omega$ is
equal to the potential function $\omega_{D,\beta}$ of the
condenser $(D,\beta)$ and therefore
\begin{equation}\label{2.14}
{\rm cap}(D,\beta)=\int_{D\setminus \beta} |\nabla\omega|^n\;dm.
\end{equation}
Finally,  (\ref{2.11}) follows from (\ref{2.12}), (\ref{2.14}) and
Lemma \ref{Lemma-6}. %
\end{proof}

Every continuum $K$ such that $\mu_D(x,y)={\rm cap}(D,K)$ for some
points $x,y\in K$ will be called $\mu_D$-minimizer. 
Simple examples show that a $\mu_D$-minimizer may not exist for
some domains $D$ and some pairs of points $x$ and $y$ and, if
exist, it is not unique, in general. We conjecture that every
$\mu_D$-minimizer $\gamma_D(x,y)$ is a smooth Jordan arc joining
$x$ and $y$.


In the last lemma of this section, we recall well-known properties
of the function $\mu_D(x,y)$, which, in particular,  show that
$\mu_D(x,y)$ is indeed a metric.

\begin{lemma}  \label{Lemma-8} %
Let $D$ be a domain in $\mathbb{R}^n$ such that $\partial D$ has
positive capacity. Then the following holds.

\begin{itemize} %
\item[(1)] %
$\mu_D(x,y)$ is a continuous function of $x$ and $y$. %
\item[(2)] $\mu_D(x,y)= 0$ if and only if $x=y$. %
\item[(3)] If $\partial D$ contains a continuum $E$, then
$\mu_D(x,y)\to \infty$ if $x\in D$ is fixed and $y\to E$. %
\item[(4)] For every triple $x,y,z$ of distinct points in $D$ the
triangle inequality holds, i.e., %
$$ 
\mu_D(x,z)\le\mu_D(x,y)+\mu_D(y,z).
$$ 
\end{itemize} %
\end{lemma} %

For the proof of properties (1) and (2) we refer to \cite{V2} and
\cite[p. 115]{F}. Property (3) follows from the monotonicity
Lemma~\ref{Lemma-5} and Lemma~7.35 in \cite{V2}. For the triangle
inequality see, for instance, \cite[p. 103]{V2}.


\section{Modulus metric and polarization} %
A geometric transformation called \emph{polarization} was
introduced by V.~Wolontis \cite{W}. Two modern approaches to this
transformation are popular now. The first one was developed by
V.~N.~Dubinin who also suggested the term ``polarization'' for
this transformation, see his book \cite{D1}, and the other
approach first appeared in \cite{S2} and then was developed in
full generality in \cite{BS}. In this paper we use polarization
with respect to spheres in $\mathbb{R}^n$, which continuously
depend on some geometric parameters. The latter approach was
inspired by papers \cite{S1} and \cite{S4}, where polarization was
used to solve P\'{o}lya-Szeg\"{o} problem on continuous
symmetrization.

\smallskip

Polarization of a set $E\subset \overline{\mathbb{R}}^n$ with
respect to a sphere $S(x_0,r)$ can be defined as follows. Given
$x\in \mathbb{R}^n\setminus \{x_0\}$, by $x^*$ we denote the point
in $\mathbb{R}^n$ symmetric to $x$ with
respect to $S(x_0,r)$, i.e. %
$$  %
x^*=x_0+r^2\frac{x-x_0}{|x-x_0|^2}. %
$$
The points $x_0$ and $x=\infty$ are considered symmetric to each
other with respect to every sphere centered at $x_0$. Let
$E^*=\{x\in \overline{\mathbb{R}}^n:\, x^*\in E\}$. Thus, $E^*$
consists of all points that are symmetric to the points of $E$
with respect to $S(x_0,r)$. In other words, $E^*$ is a reflection
of $E$ with respect to $S(x_0,r)$.

\smallskip

\noindent %
\textbf{Definition 1.} Let $E$ and $(D,K)$ be a set and a
condenser in $\mathbb{R}^n$, respectively. %
\begin{itemize} %
\item[\textbf{(a)}] Polarization of  $E$ with respect to $S(x_0,r)$ is defined as %
\begin{equation}  \label{3.1} %
E_p=\left((E\cup E^*)\cap
\overline{B_r(x_0)}\right)\cup\left((E\cap E^*)\setminus
\overline{B_r(x_0)}\right).
\end{equation} %
\item[\textbf{(b)}] Polarization of a condenser $(D,K)$ with
respect to $S(x_0,r)$ is defined as  $(D_p,E_p)$.
\end{itemize} %

It is well-known  that $E_p$ defined by (\ref{3.1}) is a compact
set if $E$ is compact and that $D_p$ is an open set if $D$ is
open. On the other side, polarization does not preserve
connectivity. Simple examples, well known to the experts, show
that there are simply connected domains $D$ the polarization of
which consists of infinitely many connected components and some of
these connected components are infinitely connected domains. Thus,
the polarization $(D_p,K_p)$ of a condenser $(D,K)$ is not, in
general, a condenser as it was defined in Section~2.
However,
everywhere below, we polarize condensers $(D,K)$ with respect to
the spheres $S(x_0,r)$ such that $B(x_0,r)\subset D$. In this
case, $D_p=D$ and the resulting pair $(D_p,K_p)=(D,K_p)$ is again
a condenser in the sense of our definition in Section~2. The
following theorem describes the
effect of polarization on the capacity of a condenser. %

\begin{theorem}[\cite{D1}]  \label{Theorem-3}%
Let $(D,K)$ be a condenser in $\mathbb{R}^n$ and $(D_p,K_p)$ be
the polarization of $(D,K)$ with respect to a sphere $S(x_0,r)$.
Suppose further that $D_p$ is connected. Then $(D_p,K_p)$ is a
condenser and %
\begin{equation} \label{3.2} %
{\mbox{cap}}(D_p,K_p)\le {\mbox{cap}}(D,K).
\end{equation}  %
\end{theorem} %

\smallskip

\noindent %
\textbf{Remark 1.} In dimension  $n=2$, the cases when equality
occurs in (\ref{3.2}) were discussed under a variety of
assumptions in \cite{D1}, \cite{S1}, and \cite{BP1}.  Also, in
dimensions $n\ge 3$, the cases of equality in polarization
inequalities for the Newtonian capacity were discussed in
\cite{S1} and \cite{BP1}. In dimensions $n\ge 3$, the question on
the cases when equality sign occurs in (\ref{3.2}), i.e. in
polarization inequality for the conformal capacity,  remains open.
Resolving this question would be an important advance in the
theory of symmetrization that also will lead to simpler proofs of
some of our results presented below.

\medskip

\noindent %
\textbf{Remark 2.}  If $K_p$ is the polarization of a connected
compact set $K$ with respect to a sphere $S(x_0,r)$, then $K_p$ is
compact but not necessarily connected. But one can easily see that
the restriction $\widetilde{K}_p=K_p\cap \overline{B(x_0,r)}$ of
$K_p$ to the closed ball $\overline{B(x_0,r)}$ is always compact
and connected.

\smallskip

Combining Theorem~3 with properties of condenser capacity
discussed in Section~2, we obtain new useful properties of the
$\mu_D$-metric presented in Lemma~\ref{Lemma-9} and
Lemma~\ref{Lemma-10} below.

\begin{lemma}  \label{Lemma-9} %
Let $D\subset \mathbb{R}^n$ be a domain  such that $\partial D$
has positive capacity. %
\begin{itemize} %
\item[(1)] Suppose that $\overline{B(x_0,r)}\subset D$. Then for
every pair of points $x$ and $y$ in $\overline{B(x_0,r)}$ there is
a $\mu_D$-minimizer $\gamma_{\mu_D}(x,y)$ which lies in
$\overline{B(x_0,r)}$.%
\item[(2)] For $x_0\in D$, let $R_0={\mbox{dist}}(x_0,\partial
D)$. Let $x_1\in B(x_0,R_0)$ and  $R_1=|x_1-x_0|$. Then, for every
$x\in \mathbb{S}$, the $\mu_D$-distance $\mu_D(x_1,x_0+tx)$,
considered as a function of $t$, is non-decreasing on $R_1\le
t<R_0$. %
\item[(3)] If $B_{\mu_D}(x_0,s)\subset B(x_0,R_0)$, then
$B_{\mu_D}(x_0,s)$ is starlike with respect to $x_0$.
\end{itemize}
\end{lemma}  %

\begin{proof} %
(1) Let $x,y\in \overline{B(x_0,r)}$ and let $\gamma_k$, $k\in
\mathbb{N}$, be a sequence of continua in $K_{xy}$ such that %
\begin{equation}  \label{3.3} %
{\mbox{cap}}\,(D,\gamma_k)\to \mu_D(x,y), \quad {\mbox{ as $k\to
\infty$.}}
\end{equation} %
Let $\gamma_k^p$ denote the polarization of $\gamma_k$ with
respect to the sphere $S(x_0,r)$ and let
$\widetilde{\gamma}_k=\gamma_k^p\cap \overline{B(x_0,r)}$. As we
mentioned above in Remark~2, $\widetilde{\gamma}_k$ is a connected
compact set in $\overline{B(x_0,r)}$ and $x,y\in
\widetilde{\gamma}_k$. Hence, $\widetilde{\gamma}_k\in K_{xy}$.
Now, it
follows from Theorem~3 and Lemma~5 that %
\begin{equation}  \label{3.4} %
\mu_D(x,y)\le {\mbox{cap}}\,(D,\widetilde{\gamma}_k)\le
{\mbox{cap}}\,(D,\gamma_k^p)\le {\mbox{cap}}\,(D,\gamma_k), \quad
{\mbox{for all $k\in \mathbb{N}$.}}
\end{equation} %
Taking the limit in (\ref{3.4}) and taking into account
(\ref{3.3}), we conclude that %
$$ %
{\mbox{cap}}\,(D,\widetilde{\gamma}_k)\to \mu_D(x,y), \quad
{\mbox{as $k\to \infty$.}}
$$ %
Now, an existence of   the required $\mu_D$-minimizer
$\gamma_{\mu_D}(x,y)\subset \overline{B(x_0,r)}$ follows from
Lemma~\ref{Lemma-7}.

\smallskip

(2) For $x\in \mathbb{S}$ and $t_1$, $t_2$ such that
$R_1<t_1<t_2<R_0$,  let $y_1=x_0+t_1x$, $y_2=x_0+t_2x$. Suppose
that $\gamma_2=\gamma_{\mu_D}(x_1,y_2)$ is a $\mu_D$-minimizer for
the points $x_1$, $y_2$. Let $\gamma_{2,p}$ denote the
polarization of $\gamma_2$ with respect to the sphere $S(x_0,r)$
with $r=\sqrt{t_1t_2}$ and let
$\widetilde{\gamma}_2=\gamma_{2,p}\cap \overline{B(x_0,r)}$. Since
the points $y_1$ and $y_2$ are symmetric with respect to
$S(x_0,r)$ it follows that $y_1\in \widetilde{\gamma}_2$. Hence,
same argument as in part (1) of this proof shows that
$\widetilde{\gamma}_2\in K_{x_1y_1}$. Therefore, applying
Theorem~3 and Lemma~5 as above,  we conclude that %
$$ %
\mu_D(x_1,y_2)={\mbox{cap}}\,(D,\gamma_{\mu_D}(x_1,y_2))\ge
{\mbox{cap}}\,(D,\gamma_{2,p})\ge
{\mbox{cap}}\,(D,\widetilde{\gamma}_2)\ge \mu_D(x_1,y_1),
$$ %
which proves the required monotonicity property.

\smallskip

(3) Now, if $B_{\mu_D}(x_0,s)\subset B(x_0,R_0)$ and $y \in
B_{\mu_D}(x_0,s)$, then $\mu_D(x_0,x_0+t(y-x_0))\le
\mu_D(x_0,y)\le s$ for all $t$, $0\le t\le 1$, by the monotonicity
property proved above. Hence, $x_0+t(y-x_0)\in B_{\mu_D}(x_0,s)$
for $0\le t\le 1$, which proves that $B_{\mu_D}(x_0,s)$ is
starlike with respect
to $x_0$. \end{proof} %


\begin{lemma}\label{Lemma-10}
Let $D$ be a domain in $\IR^{n}$  such that $\partial D$ has
positive capacity. For $x_{0}\in D$, let
$R_0={\mbox{dist}}\,(x_0,\partial D)$.  Then the function
$\mu(x)=\mu_{D}(x_{0},x)$ does not have relative extrema in the
ball $B(x_0,R_0)$ except for the absolute minimum at $x_0$.
\end{lemma} %

\begin{proof} %
(1) Suppose that there exist $x_*\not=x_0$ and $r>0$ such that
$2r<d=|x_*-x_0|$,  $\overline{B(x_*,2r)}\subset B(x_0,R_0)$, and
$\mu(x_*)\le\mu(x)$ for all $x\in \overline{B(x_*,2r)}$. By part
(1) of  Lemma~\ref{Lemma-9}, there is a continuum
$\gamma(x_*)\subset
\overline{B(x_0,d)}$ such that%
\begin{equation}  \label{3.5} %
\mu(x_*)={\mbox{cap}}\,(D,\gamma(x_*)).
\end{equation}  %
Since $\gamma(x_*)$ is closed and connected there are closed and
connected sets  $\gamma_1\subset \gamma(x_*)$ and $\gamma_2\subset
\gamma(x_*)$ satisfying the following conditions: (a)
$\gamma_1\subset \overline{B(x_0,d)}\setminus B(x_*,2r)$ and
contains the point $x_0$ and some point $x_1\in S(x_*,2r)\cap
\overline{B(x_0,d)}$, (b)   $\gamma_2\subset \overline{B(x_*,r)}$
and contains the point $x_*$ and some point $x_2\in S(x_*,r)$.

Conditions (a) and (b) show that the continua $\gamma_*$,
$\gamma_1$ and $\gamma_2$ satisfy assumptions of
Lemma~\ref{Lemma-5} concerning the cases of equality in this lemma
and therefore %
\begin{equation}  \label{3.6} %
{\mbox{cap}}\,(D,\gamma_1)<{\mbox{cap}}\,(D,\gamma_1\cup
\gamma_2)\le {\mbox{cap}}\,(D,\gamma(x_*)).
\end{equation} %
Since $\gamma_1\in K_{x_0x_1}$ we have $\mu(x_1)\le
{\mbox{cap}}\,(D,\gamma_1)$. Since $x_1\in S(x_*,2r)$ the latter
inequality combined with relations (\ref{3.5}) and (\ref{3.6})
contradicts our assumption that $\mu(x_*)\le\mu(x)$ for all $x\in
\overline{B(x_*,2r)}$. Therefore, $\mu(x)$ can not have relative
minimum in $B(x_0,R_0)$ except for the absolute minimum at $x_0$.

\smallskip

(2) Suppose that there exist $x^*$ and $r>0$ such that
$\overline{B(x^{*},2r)}\subset B(x_0,R_0)$ and
$\mu(x)\le\mu(x^{*})$ for all $x\in B(x^{*},r)$. Let $H$ be the
hyperplane passing through $x_{0}$ and orthogonal to $x^{*}-x_0$.
Below we use the following notations: $d=|x^*-x_0|$,
$\rho=|x-x_0|$, and $L=\sqrt{\rho^2+d^2}+r$; see Figure~1, which
illustrates notations used in the proofs of this section.

\begin{tikzpicture}[scale=0.8]

\draw [decoration={markings,mark=at position 1 with
{\arrow[scale=2]{>}}}, postaction={decorate},shorten >=0.4pt]
(-7,0) -- (7,0);

\draw [decoration={markings,mark=at position 1 with
{\arrow[scale=2]{>}}}, postaction={decorate},shorten >=0.4pt]
(0,-6.4) -- (0,6.4);

\draw [black,thick] (0,0) circle [radius=6];;

\node [left] at (0,-0.2) {$x_{0}$};
\node at (0,0) [circle,fill,inner sep=1.5pt]{}; 

\draw [decoration={markings,mark=at position 1 with
{\arrow[scale=2]{>}}},
postaction={decorate},shorten >=0.4pt] (0,0) --(-1.2,-5.9); 

\node [below] at (-1.2,-3.8) {$R_{0}$};

\node [below] at (1.5,0) {$x$};
\node at (1.5,0) [circle,fill,inner sep=1.5pt]{}; 

\node [below] at (0.7,0) {$\rho$};

\draw [decoration={markings,mark=at position 1 with
{\arrow[scale=2]{>}}},
postaction={decorate},shorten >=0.4pt] (1.5,0) --(-0.85,3.22); 

\draw [black,thick] (1.5,0) circle [radius=4];;

\node [left] at (0.1,1.9) {$x^{*}$};
\node at (0,2.05) [circle,fill,inner sep=1.5pt]{}; 

\draw [black,thick] (0,2.0) circle [radius=1.5];; \draw
[black,thick] (0,2.0) circle [radius=3];;

\node [left] at (0.1,1.15) {$d$};

\node [left] at (-0.3,2.5) {$r$};

\node [left] at (2.6,2.3) {$2r$};

\draw [black,thick] [domain=90:125] plot ({0.3*cos(\x)},
{2.05+0.3*sin(\x)}); \node [above] at (-0.14,2.25) {$\alpha$};

\draw [decoration={markings,mark=at position 1 with
{\arrow[scale=2]{>}}},
postaction={decorate},shorten >=0.4pt] (0,2.05) --(3,2.05); 

\draw [decoration={markings,mark=at position 1 with
{\arrow[scale=2]{>}}},
postaction={decorate},shorten >=0.4pt] (1.5,0) --(5.5,0); 

\node [above] at (4,0) {$L$};

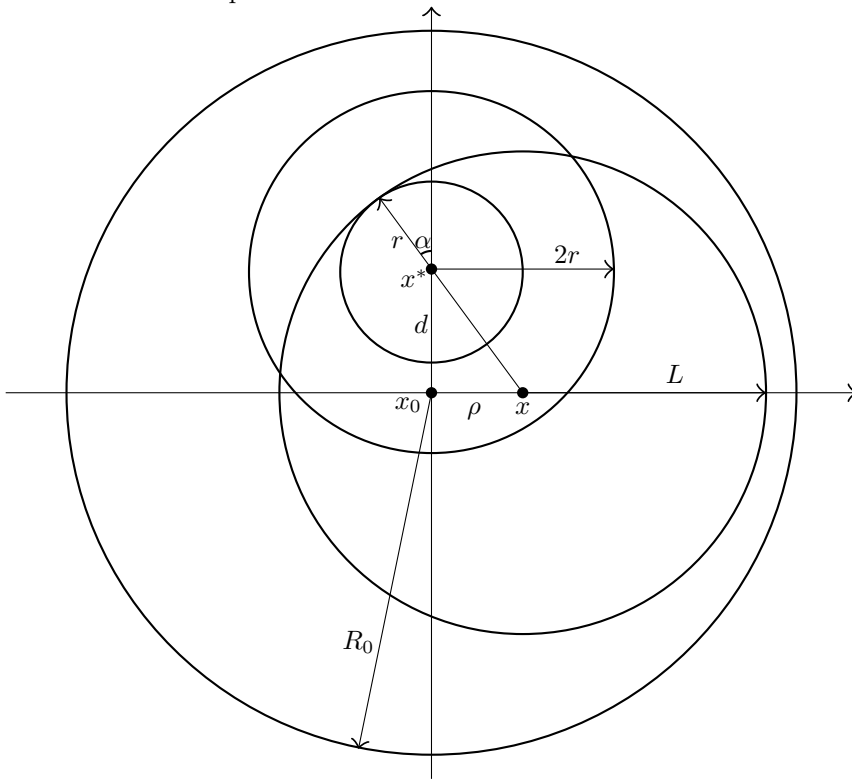
\node [below] at (-0.1,-6.7) {Figure 1. Spheres of Lemma 10.};

\end{tikzpicture}

An elementary geometric calculation shows that if $x\in H$ is such that %
\begin{equation} \label{3.7} %
\rho=|x-x_0|<\frac{r(2d+r)}{2(d+r)},
\end{equation} %
then %
\begin{equation} \label{3.8} %
 x_0\in B(x,L) \quad{\mbox{and}} \quad \overline{B(x^*,r)}\subset \overline{B(x,L)}\subset B(x_0,R_0).
\end{equation} %
Let %
\begin{equation} \label{3.9} %
y_t=x^*+t\frac{x^*-x}{|x^*-x|}, \quad 0\le t\le r.
\end{equation} %
Since conditions (\ref{3.8}) are satisfied it follows from Lemma~8
that $\mu_D(x_0,y_t)$ considered as a function of $t$ is a
non-decreasing function on $0\le t\le r$. Furthermore, since
$\mu_D(x_0,x)$ attained its relative maximum  at $x=x^*$ it
follows that $\mu_D(x_0,x)$ is constant on every radial segment of
the ball $B(x^*,r)$ of the form (\ref{3.9}) if $x$ satisfies
condition (\ref{3.7}).  Let $\Phi=\Phi(x^*,x_0,r,\alpha)$ denote
the spherical cone, which has a vertex at $x^*$, radius $r$, and
forms a central angle of opening $\alpha=\arctan
\frac{r(2d+r)}{2d(d+r)}$ with the segment
$\{y=x^*+t\frac{x^*-x_0}{|x^*-x_0|}:\, 0\le t\le r\}$. The latter
segment is a radius of the ball $B(x^*,r)$. Since every end point
$x\in S(x^*,r)$ of the radial segment from $x^*$ to $x$, which is
in the spherical cone $\Phi$, satisfies condition (\ref{3.7}) it
follows that $\mu_D(x_0,x)$ is constant on $\Phi$. Obviously,
$\Phi$ has interior points and the latter conclusion contradicts
the fact established in part (1) of this proof that $\mu(x)$ can
not have relative minimum in $B(x_0,R_0)\setminus\{x_0\}$. Thus,
our assumption was wrong and $\mu_D(x_0,x)$ does not have relative
maxima in $B(x_0,R_0)$.
\end{proof}

\smallskip

\noindent %
\textbf{Remark 3.} We conjecture that the modulus metric
$\mu(x)=\mu_D(x_0,x)$ considered as a function of $x\in D$ can not
have relative minima or relative maxima  at any point $x\in D$,
$x\not=x_0$. We want to stress here that our proof of
Lemma~\ref{Lemma-10} is based on polarization and therefore it can
not be applied to all points $x\in D$ because polarization changes
the domain $D$, in general.

\smallskip

 In the proof below, we
will use the following notations. Let $x_0$ and $x$ be points in
$D$, let $\mu_D(x_0,x)=\mu$, and let $0<|x-x_0|=r<R$, where
$R={\mbox{dist}}(x_0,\partial D)$. Also, let %
\begin{equation} \label{3.10} %
\alpha_0=\arctan\frac{\sqrt{R^2-r^2}}{r}
\end{equation} %
and %
\begin{equation} \label{3.11} %
\rho_{ext}(\alpha)=R-r\sec \alpha  \quad {\mbox{and}} \quad
\rho_{int}(\alpha)=\frac{r(R\cos\alpha-r)}{R-r\cos\alpha}.
\end{equation} %

 For $0<\alpha<\alpha_0$, by $\Phi_{ext}(\alpha)$ we denote the
spherical cone with the vertex at $x$ and  radius $\rho_{ext}$
that forms the central angle of opening $\alpha$ with the vector
$v=x-x_0$. Similarly, by $\Phi_{int}(\alpha)$ we denote the
spherical cone with the vertex at $x$ and  radius $\rho_{int}$
that forms the central angle of opening $\alpha$ with the vector
$-v_1=x_0-x$. We will call $\Phi_{ext}(\alpha)$ and
$\Phi_{int}(\alpha)$ the exterior cone and the interior cone of
$S_{\mu_D}(x_0,\mu)$, respectively. Now we are ready to prove the
cone property of the $\mu_D$-spheres stated in Theorem~2 in the
Introduction.
\smallskip

\noindent %
\textbf{Proof of Theorem~2.} %
\textbf{(1)} We claim that, for every $\alpha$,
$0<\alpha<\alpha_0$, $\Phi_{ext}(\alpha)\subset D\setminus
\overline{B_{\mu_D}(x_0,\mu)}$. Let $H_{1/2}$ be a hyperplane
passing through the point $x_{1/2}=\frac{1}{2}(x+x_0)$ and
orthogonal to the vector $v=x-x_0$. Let $x_{1/2}\in H_{1/2}$ be
such that $|x_{1/2}-x_0|=R/2$. Then the angle formed by the
vectors $v$ and $v_{1/2}=x_{1/2}-x_0$ equals $\alpha_0$; this is
how the value of $\alpha_0$ in the formula (\ref{3.10}) was
calculated.

Suppose now that for some angle $\alpha<\alpha_0$ there is a point
$y\in \Phi_{ext}(\alpha)$ that is in
$\overline{B_{\mu_D}(x,\mu)}$. By Lemma~9, the function
$\mu_D(x_0,x)$ is continuous and can not have relative minimum.
Therefore, there is a point $y^*\in \Phi_{ext}(\alpha)$ such that
$\mu_D(x_0,y^*)<\mu$.

Let $l$ be the line passing through the points $y^*$ and $x$ and
let $l$ intersects $H_{1/2}$ at the point $x^*$.  Then the angle
$\alpha^*$ formed by the vectors $v^*=x-x^*$ and $v=x-x_0$ is less
than $\alpha$; i.e., $0\le \alpha^*< \alpha<\alpha_0$. Consider
the sphere $S(x^*,\rho)$ with $\rho=\sqrt{|y^*-x^*||x-x^*|}$. The
points $x$ and $y^*$ are symmetric with respect to $S(x^*,\rho)$.
Notice also that $x_0\in B(x^*,\rho)$ and $B(x^*,\rho)\subset
B(x_0,R_0)$. The latter inclusion follows, after simple
calculations, from our definition of the radius
$\rho_{ext}(\alpha)$ defined in (\ref{3.11}).

Let $\gamma_{\mu_D}(x_0,y^*)$ be a $\mu_D$-minimizer for the
points $x_0$ and $y^*$ and let $\gamma_p\subset
\overline{B(x^*,\rho)}$ denote the connected component of the
polarization of $\gamma_{\mu_D}(x_0,y^*)$ with respect to the
sphere $S(x^*,\rho)$. Then $\gamma_p$ is a continuum in $D$
containing the points $x_0$ and $x$ and therefore $\gamma_p\subset
K_{x_0x}$. Now, using Theorem~3, we obtain the following %
$$ %
\mu_D(x_0,y^*)={\mbox{cap}}\,(D,\gamma_{\mu_D}(x_0,y^*))\ge
{\mbox{cap}}\,(D,\gamma_p)\ge \mu_D(x_0,x)=\mu.
$$ %
The latter inequalities contradicts our assumption that
$\mu_D(x_0,y^*)<\mu$, which proves our claim on the exterior
spherical cone.

\smallskip

\textbf{(2)} Now, we prove that and $\Phi_{int}(\alpha)\subset
B_{\mu_D}(x_0,\mu)$. The proof is similar to the proof given in
part \textbf{(1)}. We assume that there is a point $y$ in
$\Phi_{int}(\alpha)$ that is not in $B_{\mu_D}(x_0,\mu)$. Since,
by Lemma~11, $\mu_D(x_0,x)$ does not have relative maxima in
$B(x_0,R_0)$ it follows that there is a point $y^*$ in
$\Phi_{int}(\alpha)$ such that $\mu_D(x_0,y^*)>\mu$. As in part
\textbf{(1)}, we consider the line $l$ passing though the points
$y^*$ and $x$. Let $\alpha^*$ denote the angle formed by the
vectors $v^*=y^-x$ and $-v=x_0-x$. Then $0\le
\alpha^*<\alpha<\alpha_0$. Let $x^*$ be the point on $l$ such that
$|x^*-x|=R_0-|x^*-x_0|$. Consider the sphere $S(x^*,\rho)$ with
the radius $\rho=\sqrt{|x-x^*||y^*-x^*|}$. Then the points $y^*$
and $x$ are symmetric with respect to $S(x^*,\rho)$. Furthermore,
$x_0\in B(x^*,\rho)\subset B(x_0,R_0)$.

Let $\gamma_{\mu_D}(x_0,x)$ be a $\mu_D$-minimizer for points
$x_0$ and $x$ and let $\gamma_p\subset \overline{B(x^*,\rho)}$
denote the connected component of the polarization of
$\gamma_{\mu_D}(x_0,x)$ with respect to the sphere $S(x^*,\rho)$.
Then $\gamma_p$ is a continuum in $D$ containing the points $x_0$
and $y^*$ and therefore $\gamma_p\subset
K_{x_0y^*}$. Now, using Theorem~3, we obtain the following %
$$ %
\mu=\mu_D(x_0,x)={\mbox{cap}}\,(D,\gamma_{\mu_D}(x_0,x))\ge
{\mbox{cap}}\,(D,\gamma_p)\ge \mu_D(x_0,y^*).
$$ %
The latter inequalities contradicts our assumption that
$\mu_D(x_0,y^*)>\mu$, which proves our claim on the interior
spherical cone. %
The proof of Theorem~2 is complete.  \hfill $\Box$

\smallskip

It follows from Theorem~2 and its proof above that stronger
versions of statements of Lemma~\ref{Lemma-9} hold true. We
present these stronger versions in the following corollary.

\begin{corollary}  \label{Corollary-1} %
Under the assumptions of Theorem 2, the following statements hold true.%
\begin{itemize}  %
\item[(1)] Let $L(\alpha_0)$ denote the spherical lune which is
the intersection of all balls $B(y,R_0/2)$ having centers at $y\in
H_{1/2}$ such that $|y-x_0|=R_0/2$. Then there is a
$\mu_D$-minimizer $\gamma_{\mu_D}(x_0,x)$ contained in
$L(\alpha_0)$. %
 \item[(2)]  Let $y\in \mathbb{S}$. The function
$\mu_D(t)=\mu_D(x_0,x_0+ty)$ is strictly increasing for $0\le t<R_0$. %
In particular, $S_{\mu_D}(x_0,\tau)$ can not contain intervals of
a line passing through the point $x_0$.%
 \item[(3)] If
$\overline{B_{\mu_D}(x_0,\tau)}\subset B(x_0,R_0)$, then
$B_{\mu_D}(x_0,\tau)$ is strictly starlike, i.e. every ray $l_+$
in $\mathbb{R}^n$ with the initial point at $x_0$ intersect the
$\mu_D$-sphere  $S_{\mu_D}(x_0,\tau)$ at one point. %
\end{itemize}  %
\end{corollary} %

\begin{proof} %
Part (1) follows from the fact used in the proof of Theorem~2 that
every ball $B(y,R_0/2)$ with the center $y\in H_{1/2}$ such that
$|y-x_0|=R_0/2$ contains a $\mu_D$-minimizer. Parts (2) and (3)
follow immediately from  the cone properties of the
$\mu_D$-spheres. \end{proof} %

\section{Infinitesimally small spheres and conformality}

Everyone who studied Complex Analysis remembers that conformal
mapping transforms small circles to ``\emph{infinitesimally small
circles}''. However, it was not easy for us to find a precise
definition of this term, especially in the $n$-dimensional
setting,  in the accessible literature. For our purposes, we adapt
the definition introduced by Yu. G. Reshetnyak \cite{R1}.

\smallskip

\noindent %
\textbf{Definition 2.}  Let $D$ be a domain in $\mathbb{R}^n$ and
$x_0\in D$.
\begin{itemize} %
\item[(1)] %
A parameterized family ${\mathcal{U}}(x_0)=\{U_t(x_0):\, 0<t\le
t_0\}$, of neighborhoods $U_t(x_0)\subset D$ of $x_0$ is called
\emph{almost spherical} if the
following holds: 
\begin{itemize} %
\item[(a)] $U_t(x_0)\subset U_{t_0}(x_0)$ for all $t$ and there is
a homeomorphism $\varphi$ from $U_{t_0}(x_0)$ to $\mathbb{R}^n$
such that $\varphi(U_t(x_0))=S(\bar{0},t)$ for $0< t\le
t_0$. %
\item[(b)] $\max_{x\in \partial U_t(x_0)} |x-x_0|/\min_{x\in
\partial U_t(x_0)}
|x-x_0|\to 1$ as $t\to 0$.  %
\end{itemize} %
\item[(2)] If ${\mathcal{U}}(x_0)=\{U_t(x_0):\, 0<t\le t_0\}$ is
an almost spherical family of neighborhoods $U_t(x_0)$ in $D$,
then the family $\Sigma(x_0)=\{\partial U_t(x_0):\,0<t\le t_0\}$
consisting of the boundary surfaces of $U_t(x_0)$ will be called
an \emph{infinitesimally small sphere} centered at $x_0$.
\end{itemize} %

\medskip

With this terminology, the main result of Reshetnyak's paper
\cite{R1} can be stated in the following form. %
\begin{theorem}[\cite{R1}] %
Let $D$ be a domain in $\mathbb{R}^n$ and let
$\mathcal{S}=\{\Sigma(x_0):\,x_0\in D\}$ be a collection of
infinitesimally small spheres $\Sigma(x_0)$ centered at $x_0$ such
that one such sphere is assigned to each point $x_0\in D$. Then a
homeomorphism $f$ from $D$ onto a domain $\Omega\in
\mathbb{R}^n$ is conformal if and only if for every $x_0\in D$ %
the image $f(\Sigma(x_0))$ is an infinitesimally small sphere in
$\Omega$ centered at $f(x_0)$.
\end{theorem}  %

In view of this Reshetnyak's theorem, to prove Theorem~1 we have
to show that for every domain $D\subset \mathbb{R}^n$ and every
point $x_0\in D$ an appropriate truncation
$\Sigma'(x_0)=\{S_{\mu_D}(x_0,t):\, 0<t\le t_0(x_0)\}$ of the
family of level sets of the modulus metric $\mu_D$ is an
infinitesimally small sphere in $D$ centered at $x_0$. This will
be established in Lemma~\ref{Lemma-12} below.
To prove this lemma, 
we will use polarization with respect to appropriate
spheres. An existence of such spheres follows from our next lemma
that can be seen as an exercise in elementary geometry.

\begin{lemma}[Lemma on $3$ spheres] \label{Lemma-11} %
Let $S_1$ and $S_2$ be two
concentric spheres centered at $x_0\in \mathbb{R}^n$ of radii
$R_1=R$ and $R_2=kR$, respectively, with $R>0$ and $0<k<1$. Then
for every pair of points $x_1\in S_1$ and $x_2\in S_2$ there is a
sphere $S_3$ of radius $R_3=R_3(x_1,x_2)$ such
that: %
\begin{itemize} %
\item[(1)] $x_1$ and $x_2$ are symmetric with respect to $S_3$. %
\item[(2)] $x_0$ and $x_2$ belong to the closed ball bounded by $S_3$. %
\item[(3)] $\frac{k}{\sqrt{1+k^2}}R\le R_3(x_1,x_2)\le
\frac{k}{1-k}R$.
\end{itemize}
\end{lemma} %
\begin{proof}%
Let $x_1\in S_1$ and $x_2\in S_2$. Using translation and scaling,
if necessary, we may assume without loss of generality that $S_1$
is the sphere of radius $1$ centered at $x_0=\bar{0}$, then $S_2$
is the sphere of radius $k$ centered at $x_0=\bar{0}$.
Furthermore, using rotation and reflection, again if necessary, we
may assume
 that $x_1$ and $x_2$ belong
to a two-dimensional plane $P$, embedded in $\mathbb{R}^n$, 
and that in the plane $P$ the points $x_1$ and $x_2$ have the
following two-dimensional coordinates:
$x_1=(-\sqrt{1-k^2\sin^2\theta},k\sin \theta)$ and
$x_2=(k\cos\theta,k\sin\theta)$, $0\le \theta\le \pi$. Thus, under
these assumptions,  the points $x_1$ and $x_2$ lie on the same
horizontal line $L\subset P$. See Figure~2, which illustrates
notations used in this proof.

\begin{tikzpicture}[scale=0.8]
\draw [decoration={markings,mark=at position 1 with
{\arrow[scale=2]{>}}}, postaction={decorate},shorten >=0.4pt]
(-7,0) -- (7,0);

\draw [decoration={markings,mark=at position 1 with
{\arrow[scale=2]{>}}}, postaction={decorate},shorten >=0.4pt]
(0,-5) -- (0,6);

\draw [black,thick] (0,0) circle [radius=2.5];; \draw
[black,thick] (0,0) circle [radius=4];;
\draw (-7,1.7) --(7,1.7); 
\node [above] at (6.8,1.7) {$L$};
\node at (0,0) [circle,fill,inner sep=1.5pt]{}; 
\node at (1.82,1.7) [circle,fill,inner sep=1.5pt]{}; 
\node at (-3.62,1.7) [circle,fill,inner sep=1.5pt]{}; 
\node [left] at (0,-0.2) {$x_{0}$}; \node [above] at (1.84,1.7)
{$x_{2}$}; \node [above] at (-3.72,1.7) {$x_{1}$}; \draw (0,0)
--(1.82,1.7); \draw (0,0) --(-3.62,1.7);
\node at (5.5,1.7) [circle,fill,inner sep=1.5pt]{}; 
\node [above] at (5.5,1.7) {$x_{c}$}; \node [above] at (3,0.9)
{$R_{3}$};

\draw  [decoration={markings,mark=at position 1 with
{\arrow[scale=2]{>}}},
postaction={decorate},shorten >=0.4pt](5.5,1.7) --(0,0); 

\draw [decoration={markings,mark=at position 1 with
{\arrow[scale=2]{>}}},
postaction={decorate},shorten >=0.4pt] (0,0) --(3.89,-0.9); 

\node [below] at (3.35,-0.7) {$R_{1}$};

\draw [decoration={markings,mark=at position 1 with
{\arrow[scale=2]{>}}},
postaction={decorate},shorten >=0.4pt] (0,0) --(-1.2,-2.2); 

\node [below] at (-1.1,-0.9) {$R_{2}$};

\draw [black,thick] [domain=0:43] plot ({0.74*cos(\x)},
{0.74*sin(\x)}); \node [above] at (0.8,0.2) {$\theta$};

\node [below] at (-1.5,2.7) {$S_{2}$}; \node [below] at
(-2.46,3.9) {$S_{1}$}; \node [below] at (3.1,-3.5) {$S_{3}$};
\node [below] at (-4.1,4) {$P$};

\draw [black,thick] [domain=160:240] plot ({5.5+5.74*cos(\x)},
{1.7+5.74*sin(\x)}); \draw [black,thick] [domain=144:146] plot
({5.5+5.74*cos(\x)}, {1.7+5.74*sin(\x)}); \draw [black,thick]
[domain=148:150] plot ({5.5+5.74*cos(\x)}, {1.7+5.74*sin(\x)});
\draw [black,thick] [domain=152:154] plot ({5.5+5.74*cos(\x)},
{1.7+5.74*sin(\x)}); \draw [black,thick] [domain=156:158] plot
({5.5+5.74*cos(\x)}, {1.7+5.74*sin(\x)});

\draw [black,thick] [domain=242:244] plot ({5.5+5.74*cos(\x)},
{1.7+5.74*sin(\x)}); \draw [black,thick] [domain=246:248] plot
({5.5+5.74*cos(\x)}, {1.7+5.74*sin(\x)}); \draw [black,thick]
[domain=250:252] plot ({5.5+5.74*cos(\x)}, {1.7+5.74*sin(\x)});
\draw [black,thick] [domain=254:256] plot ({5.5+5.74*cos(\x)},
{1.7+5.74*sin(\x)});

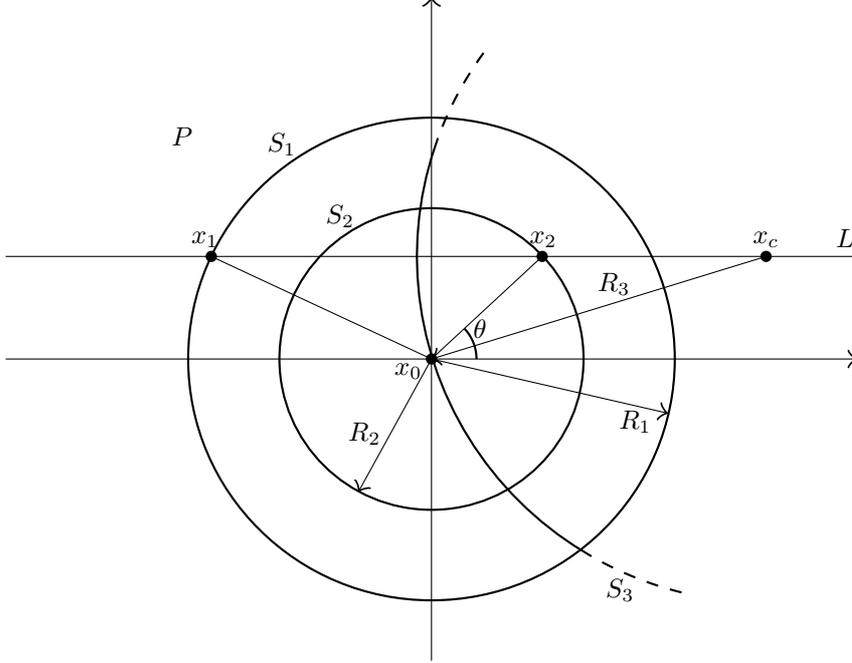
\node [below] at (-0.1,-5.3) {Figure 2. Three spheres lemma.};
\end{tikzpicture}

1) First, we assume that $\pi-{\mbox{arccot}}\,k\le\theta\le \pi$.
In this case, we define $S_3$ to be a sphere in $\mathbb{R}^n$
centered at the point $x_c\in P$  with coordinates $(0,k\sin
\theta)$ in the plane $P$ and radius %
\begin{equation} \label{4.1} %
R_3=\sqrt{-k\cos \theta\sqrt{1-k^2\sin^2\theta}}.
\end{equation}  %
The latter equation shows that the points $x_1$ and $x_2$ are
symmetric with respect to $S_3$. Furthermore, an elementary
calculation shows that  for all $\theta$,
$\pi-{\mbox{arccot}}\,k\le \theta\le \pi$,
the following equation holds:%
\begin{equation} \label{4.2} %
R_3^2=-k\cos \theta \sqrt{1-k^2\sin^2\theta}\ge
k^2\sin^2\theta=({\mbox{dist}}(x_c,x_0))^2. %
\end{equation}%
Moreover, equality occurs in (\ref{4.2}) only for
$\theta=\pi-{\mbox{arccot}}\,k$. Thus, inequality (\ref{4.2})
implies that the points $x_2$ and $x_0=\bar{0}$ belong to a closed
ball in $\mathbb{R}^n$ bounded by $S_3$. In particular, if
$\theta=\pi-{\mbox{arccot}}\,k$, then the sphere $S_3$ passes
through the origin $x_0=\overline{0}$. Therefore the sphere $S_3$
satisfies conditions $(1)$ and $(2)$ for the values of $\theta$
under consideration. In addition, it is immediate from (\ref{4.1})
that the radius $R_3$ is an
increasing function of $\theta$. Hence, %
\begin{equation} \label{4.3} %
\frac{k}{\sqrt{1+k^2}}\le R_3=R_3(x_1,x_2)\le \sqrt{k}, \quad
\quad {\mbox{when $\pi-{\mbox{arccot}}\,k\le \theta\le
\pi$}}. %
\end{equation} %

2) Now we turn to the case when $0\le \theta\le
\pi-{\mbox{arccot}}\,k$. We claim that in this case there is a
unique point $x_c\in L$, with coordinates
$(\lambda_c,k\sin\theta)$, $\lambda_c>0$, in the plane $P$, such
that the points $x_1$ and $x_2$ are symmetric with respect to the
sphere
$S_3=S_3(x_1,x_2)$ centered at $x_c$ with radius %
$ 
R_3=\sqrt{\lambda_c^2+k^2\sin^2\theta}.
$ 
Notice  that under these conditions, the sphere $S_3$ passes
through the origin $x_0=\bar{0}$.

First we introduce necessary notations. We fix $\theta$, $0\le
\theta<\pi-{\mbox{arccot}}\,k$, and consider a point $x(\rho)$ on
the line $L$, which is uniquely
determined by the following conditions: %
\begin{itemize} %
\item[(a)] $x(\rho)$ lies further to the right on $L$ than
$x_2$ and the point  $x_i=(0,k\sin\theta)$. %
\item[(b)] The distance from $x(\rho)$ to the origin
$x_0$ equals $\rho$.  %
\end{itemize} %

Let $S(\rho)$ denote the sphere centered at $x(\rho)$ such that
the points $x_1$ and $x_2$ are symmetric with respect to
$S(\rho)$. Then the radius $\tau=\tau(\rho)$ of this sphere can be
found
from the equation %
\begin{equation} \label{4.4} %
\tau^2=(\sqrt{\rho^2-k^2\sin^2\theta}-k\cos\theta)\left(\sqrt{\rho^2-k^2\sin^2\theta}+\sqrt{1-k^2\sin^2\theta}\right).
\end{equation}  %

To prove the existence part of our claim, it is enough  to show
that the equation $\tau(\rho)=\rho$ has at least one solution in
the case under consideration. The existence of such solution
follows from the continuity of function $\tau(\rho)$ given by
(\ref{4.4}) and the following ``boundary'' relations: %
\begin{enumerate} %
\item[($\alpha$)] %
If $0\le \theta\le \pi/2$ and $\rho=k$, then
$\tau(\rho)=\tau(k)=0<k=\rho$. %
\item[($\alpha'$)] %
If $\pi/2< \theta< \pi-{\mbox{arccot}}\,k$ and $\rho=k\sin\theta$,
then it follows from our argument in part (1) of this proof that
$\tau(\rho)=\tau(k\sin\theta)<k\sin\theta=\rho$. %
\item[($\beta$)] Using equation (\ref{4.4}) one can easily show
that the function $\tau^2(\rho)$ has the following asymptotic
expansion:  %
$$ %
\tau^2(\rho)=
\rho^2+(\sqrt{1-k^2\sin^2\theta}-k\cos\theta)+o(\rho), %
$$ %
where $o(\rho)\to 0$ when $\rho\to \infty$. %
\end{enumerate}  %

 Relations ($\alpha$)
and ($\beta$) show that the difference $\tau(\rho)-\rho$ changes
its sign when $\rho$ runs from $k$ to $\infty$; similarly,
relations ($\alpha'$) and ($\beta$) show that $\tau(\rho)-\rho$
changes its sign when $\rho$ runs from $\pi-{\mbox{arccot}}\,k$ to
$\infty$. Therefore, in each of these cases equation
\begin{equation}\label{4.5}  %
\rho^2=
(\sqrt{\rho^2-k^2\sin^2\theta}-k\cos\theta)(\sqrt{\rho^2-k^2\sin^2\theta}+\sqrt{1-k^2\sin^2\theta})  %
\end{equation} %
has at least one solution in the corresponding interval. In fact,
equation (\ref{4.5}) can be easily solved and its unique solution
$R_3=R_3(x_1,x_2)$ is  %
\begin{equation}  \label{4.6} %
R_3^2=\left(\frac{k^2}{\sqrt{1-k^2\sin^2\theta}-k\cos\theta}+k\cos\theta\right)^2+k^2\sin^2\theta.
\end{equation} %

 Differentiating both sides of equation (\ref{4.6}) with respect
 to $\theta$  and then simplifying the output, we obtain: %
\begin{equation} \label{4.7} %
2R_3\frac{dR_3}{d\theta}=-\frac{2k^3\sin\theta}{\sqrt{1-k^2\sin^2\theta}\left(\sqrt{1-k^2\sin^2\theta}-k\cos\theta\right)^2}<0. %
\end{equation}  %

Since the derivative $\frac{dR_3}{d\theta}$ is negative, the
radius $R_3=R_3(\theta)$, considered as a function of $\theta$,
strictly decreases from $R_3(0)=\frac{k}{1-k}$ to
$R_3(\pi-{\mbox{arccot}}\,k)=\frac{k}{\sqrt{1+k^2}}$ when $\theta$
varies from $0$ to $\pi-{\mbox{arccot}}\,k$. The latter together
with the inequality (\ref{4.2}) proves part (3) of the lemma. Now
the proof is complete. \end{proof}

\noindent %
\textbf{Remark 4.} In notations of Lemma~\ref{Lemma-11},
suppose that $x_1=(-\sqrt{1-k^2\sin^2\theta},k\sin \theta)$,
$x_2=(k\cos\theta,k\sin\theta)$ and that $0\le \theta\le \pi/2$.
In this case, the monotonicity property of the radius
$R_3(x_1,x_2)$ established in the proof of Lemma~\ref{Lemma-11}
implies the
following bounds for this radius: %
\begin{equation} \label{4.8}
\frac{k}{\sqrt{1-k^2}}R\le R_3(x_1,x_2)\le \frac{k}{1-k}R.
\end{equation} %
Thus, $R_3(x_1,x_2)\to \infty$ uniformly on $0\le \theta \le
\pi/2$ as $k\to 1$.

\begin{lemma}\label{Lemma-12}
Let $D$ be a domain in $\IR^{n}$. For every $x_{0}\in D$, the
family $\Sigma(x_0)=\{S_{\mu_D}(x_{0},t):t\in (0,\infty)\}$ of the
level surfaces of the modulus metric $\mu_D(x_0,y)$ considered as
a function of $y\in D$ has a truncation
$\Sigma'(x_0)=\{S_{\mu_D}(x_0,t):\, 0<t\le t_0(x_0)\}$ which is
an infinitesimally small sphere centered at $x_{0}$.
\end{lemma} %
\begin{proof} %
Fix $x_0\in D$. We have to show that an appropriate truncation of
$\Sigma(x_0)$ satisfies conditions (a) and (b) of part (1) of
Definition~2.

Since, by Lemma~\ref{Lemma-8}, the function $\mu_D(x_0,x)$ is
continuous and $\mu_D(x_0,x)\to 0$ as $x\to x_0$, there is $t_0>0$
such that $S_{\mu_D}(x_0,t)\subset S(x_0,R_0)$ for all $t$,
$0<t\le t_0$. Here $R_0={\mbox{dist}}\,(x_0,\partial D)$. We claim
that $\Sigma'(x_0)=\{S_{\mu_D}(x_0,t):\,0<t\le t_0\}$ is an
infinitesimally small sphere.

Indeed, consider the mapping $\varphi$ defined by %
\begin{equation}  \label{4.9} %
\varphi(x)=\frac{x}{|x|}\, \mu_D(x_0,x), \quad \quad x\in
B_{\mu_D}(x_0,t_0).
\end{equation} %

Since $\mu_D(x_0,x)$ is continuous by part (1) of
Lemma~\ref{Lemma-8} and it is strictly increasing by part (2) of
Corollary~1, it follows from (\ref{4.9} ) that $\varphi$ maps
$B_{\mu_D}(x_0,t_0)$ continuously and one-to-one onto the ball
$B(\bar{0},t_0)$ and such that
$\varphi(S_{\mu_D}(x_0,t))=S(\bar{0},t)$ for all $t$, $0<t\le
t_0$. Therefore, the family
${\mathcal{U}}(x_0)=\{B_{\mu_D}(x_0,t):\, 0<t\le t_0\}$ of
neighborhoods of $x_0$ and the mapping $\varphi$ satisfy condition
(a) of part (1) of Definition~2.

\smallskip

It remains to show that the family $\Sigma'(x_0)$ satisfies
condition (b) of part (1) of Definition~2. Suppose that this
condition is not satisfied for a sequence of the $\mu_D$-spheres
$S_i=S_{\mu_D}(x_0,t_i)$, $i\in \mathbb{N}$, where $t_i\to 0$ as
$i\to \infty$. Then there are an index $i_0\in \mathbb{N}$ and
$k$, $0<k<1$, such that for every $i\ge i_0$ there are points
$x_i,y_i\in S_{\mu_D}(x_0,t_i)$ such that %
\begin{equation} \label{4.10} %
\frac{|x_0-x_i|}{|x_0-y_i|}=k \quad {\mbox{for all $i\ge i_0$.}}
\end{equation} %
Furthermore, since, by Lemma~\ref{Lemma-10}, $\mu_D(x_0,x)$ does
not have relative maxima in $B(x_0,R_0)$ it follows that for every
$i\ge
i_0$ there is $z_i\in B(x_0,R_0)$ such that %
\begin{equation} \label{4.11} %
\frac{|x_i-z_i|}{|x_0-y_i|}<\frac{1}{i} \quad {\mbox{and}} \quad
\mu_D(x_0,z_i)>t_i.
\end{equation} %
Let $k_i=|x_0-z_i|/|x_0-y_i|$. From (\ref{4.10}) and (\ref{4.11}), we conclude that the inequalities %
\begin{equation}  \label{4.12} %
\frac{k}{2} \le k_i\le \frac{1+k}{2}<1
\end{equation} %
hold for all $i\ge i_1$, if  $i_1\ge \max
\{i_0,\frac{2}{k},\frac{2}{1-k}\}$.

It follows from Lemma~\ref{Lemma-11} that for every $i\ge i_1$
there is a sphere $S_i$ such that $z_i$ and $y_i$ are symmetric
with respect to $S_i$, the points $x_0$ and $z_i$ belong to the
closed ball $\overline{B_i}$ bounded by $S_i$ and such that the
radius $R_i$
of $S_i$ satisfies the inequalities %
\begin{equation} \label{4.13} %
R_i\le \frac{k_i}{1-k_i}|y_i-x_0|\le \frac{1+k}{1-k}|y_i-x_0|,
\end{equation}  %
where the second inequality follows from (\ref{4.12}). Since
$y_i\to x_0$ as $i\to \infty$ it follows from (\ref{4.13}) that
there is $i_* \ge i_1$ such that $R_{i_*}<R_0/2$. Since $x_0\in
\overline{B_{i_*}}$, the latter implies that
$\overline{B_{i_*}}\subset B(x_0,R_0)$.

Let $\gamma_{\mu_D}(x_0,y_{i_*})$ be a $\mu_D$-minimizer for the
points $x_0$ and $y_{i_*}$ and let $\gamma_p\subset
\overline{B_{i_*}}$ denote the connected component of the
polarization of $\gamma_{\mu_D}(x_0,y_{i_*})$ with respect to the
sphere $S_{i_*}$. Then $\gamma_p$ is a continuum in $D$ containing
the points $x_0$ and $z_i$ and therefore $\gamma_p\subset
K_{x_0z_i}$. Now, using Theorem~3 and (\ref{4.11}), we obtain the following %
$$ %
\mu_D(x_0,y_{i_*})={\mbox{cap}}\,(D,\gamma_{\mu_D}(x_0,y_{i_*}))\ge
{\mbox{cap}}\,(D,\gamma_p)\ge \mu_D(x_0,z_i)>t_{i_*}.
$$ %
The latter inequalities contradicts our assumption that
$\mu_D(x_0,y_{i_*})=t_{i_*}$.

Thus, our assumption that there is a sequence
$S_i=S_{\mu_D}(x_0,t_i)$, $i\in \mathbb{N}$, of $\mu_D$-spheres,
which do not satisfy condition (b) of part (1) of Definition~2
leads to a contradiction. Therefore, the family $\Sigma'(x_0)$
satisfies this condition, which completes our proof of Lemma~\ref{Lemma-12}.%
\end{proof} %

 \smallskip


Now we are ready to prove our main result.

\smallskip

\noindent%
\textbf{Proof of Theorem 1.} Let $f:D\to \Omega$ be an isometry of
the metric space $(D,\mu_D)$ onto the metric space
$(\Omega,\mu_\Omega)$. For   $x_0\in D$, let
$R_D(x_0)={\mbox{dist}}\,(x_0,\partial D)$ and
$R_\omega(x_0)={\mbox{dist}}(f(x_0),\partial \Omega)$. To each
$x_0\in D$ we assign a family
$\Sigma'(x_0)=\{S_{\mu_D}(x_0,t):\,0<t\le t_0\}$ of
$\mu_D$-spheres $S_{\mu_D}(x_0,t)$ such that
$S_{\mu_D}(x_0,t)\subset B(x_0,R_D(x_0))$ and
$f(S_{\mu_D}(x_0,t))\subset B(f(x_0),R_\Omega(f(x_0)))$ for all
$t$, $0<t\le t_0$. It follows from Lemma~\ref{Lemma-12} that, for
each $x_0$, $\Sigma'(x_0)$ is an infinitesimally small sphere in
$D$ centered at $x_0$. Also, since $f$ is an isometry from
$(D,\mu_D)$ to $(\Omega,\mu_\Omega)$ it follows that
$f(S_{\mu_D}(x_0,t))=S_{\mu_\Omega}(f(x_0),t)$. Since
$S_{\mu_\Omega}(f(x_0),t)=f(S_{\mu_D}(x_0,t))\subset
B(f(x_0),R_\omega(f(x_0)))$ for all $t$, $0<t\le t_0$, it follows
from Lemma~\ref{Lemma-12} that
$\Sigma'_f(f(x_0))=\{f(S_{\mu_D}(x_0,t)):\,0<t\le t_0\}$ is an
infinitesimally small sphere in $\Omega$. Thus, for every point
$x_0$ in $D$ there is an infinitesimally small sphere
$\Sigma'(x_0)$ that is mapped by $f$ onto an infinitesimally small
sphere in $\Omega$ which is centered at $f(x_0)$. Therefore, by
Theorem~4, $f$ is a conformal mapping from $D$ to $\Omega$. \hfill
$\Box$

\smallskip

\noindent %
\textbf{Remark 5.} It is tempting to use the polarization
technique alone, without referencing to the rather deep
Reshetnyak's Theorem~4, to prove conformality of isometries $f$
between metric spaces $(D,\mu_D)$ and $(\Omega,\mu_\Omega)$. At a
first glance it looks possible since, if the image $f(S(x_0,r))$
of a sphere $S(x_0,r)\subset D$ is not a round sphere, then it is
squeezed between two spheres $S_1=S(f(x_0),r)$ and
$S_2=S(f(x_0),R)$ such that $0<r/R=k<1$. Then, by
Lemma~\ref{Lemma-11}, we may find the third sphere $S_3$ and then
use polarization with respect to $S_3$ as in the proof of
Lemma~\ref{Lemma-12} to get a contradiction to the assumption of
non-roundness of $f(S(x_0,r))$. The only obstacle for this
``proof'' is the inequality (\ref{4.8}) of Remark~4. Precisely,
this inequality shows that the radius $R_3$ of the sphere $S_3$
may grow without bounds as $k\to 1$ and therefore polarization
with respect to $S_3$ will eventually destroy the domain $D$ if it
is not the whole space $\mathbb{R}^n$.

\section{Open questions and further research}
Our polarization approach is essentially geometrical and thus can
be adapted to prove similar results for some other metrics. What
is needed is a few basic properties of the metric and polarization
inequality akin to (\ref{3.2}). But polarization (or
symmetrization) alone does not provide enough information to
study, for instance, delicate properties of the $\mu_D$-spheres
and $\mu_D$-minimizers while other tools are not available. This
is why many questions about their structure remain open. Below, we
mention three of them.

\smallskip

\noindent %
 \textbf{Problem~1.} Prove that the $\mu_D$-spheres in
$D\subset \mathbb{R}^n$, $n\ge 3$, generically are smooth
topological spheres or finite collections of disjoint smooth
topological spheres.

Describe the structure of $S_{\mu_D}(x_0,r)$ near its critical
points; i.e. near the points $x\in D$, where $S_{\mu_D}(x_0,r)$ is
not smooth.

\smallskip

\noindent %
 \textbf{Problem~2.} Prove that  every $\mu_D$-minimizer
$\gamma_{\mu_D}(x,y)$  in $D\subset \mathbb{R}^n$, $n\ge 3$, is a
smooth Jordan arc. Currently it is not known if
$\gamma_{\mu_D}(x,y)$ is an irreducible continuum or even whether
or not $\gamma_{\mu_D}(x,y)$ may have interior points.

\smallskip

To state our next problem we need some terminology. If $D$ is a
domain in $\overline{\mathbb{R}}^n$, $x,y\in D$,  and  $\gamma\in
K_{xy}$ is such that $\mu_D(x,y)={\mbox{cap}}\,(D,\gamma)$, then
we say that $\gamma$ is a $\mu_D$-minimizer with endpoints $x,y$.
We say that a family $\Gamma=\{\gamma\}$ of $\mu_D$-minimizers
$\gamma\subset D$ foliates $D$ if:
 (a) $\cup_{\gamma\in \Gamma} \gamma=D$ and (b) if
 $\gamma_1,\gamma_2\in \Gamma$ and there is $x\in \gamma_1\cap
 \gamma_2$, which is not an endpoint for at least one of these $\mu_D$-minimizers, then either $\gamma_1\subset \gamma_2$ or
 $\gamma_2\subset \gamma_1$.

\smallskip

\noindent %
 \textbf{Problem~3.} Let $D$ be a domain in $\mathbb{R}^n$, $n\ge 3$,
 supplied with the $\mu_D$-metric. Then $D$ has a family of
 $\mu_D$-minimizers foliating $D$ if and only if $D$ is a
 topological ball or $D$ is a topological spherical shell.

\smallskip

\noindent%
\textbf{Remark 6.} In all three problems stated above we assume
that $n\ge 3$. In the planar case, when $n=2$, these problems are
easier and can be 
resolved within the frame of the
Jenkins' theory on extremal partitioning, see \cite{S3}. But this
is already a topic for another paper.

\end{document}